\documentclass[11pt]{article}
\usepackage{latexsym,amssymb,amsmath,amscd,amsthm,amsxtra}
\usepackage{mathrsfs}
\usepackage{epsfig}
\usepackage{pgf,tikz}
\usepackage{graphicx}
\usepackage{graphics}
\usepackage{setspace}
\newcommand{\counte}{section}

\newtheorem{prop}{\bf Proposition}[\counte]
\newtheorem{lemma}{\bf Lemma}[\counte]
\newtheorem{theorem}{\bf Theorem}[\counte]
\newtheorem{coro}{\bf Corollary}[\counte]
\newtheorem{example}{\bf Example}[\counte]

\newtheorem{remark}{\bf Remark}[\counte]

\author{Zhao Xu-an, zhaoxa@bnu.edu.cn\\ Department of Mathematics, Beijing Normal University\\Key Laboratory
of Mathematics and Complex Systems\\ Ministry of Education,
China, Beijing 100875
\\Jin Chunhua, jinch@amss.ac.cn\\Academy of Mathematics and Systems Science, Chinese Academy of Sciences
\\China, Beijing 100190}

\title{Polynomial invariants of Weyl groups for Kac-Moody groups\thanks{The authors are supported by NSFC11171025}}

\date{}

\begin{document}

\maketitle
\begin{abstract}
In this paper, we prove that the ring of polynomial invariants of the Weyl group for an indecomposable and indefinite Kac-Moody Lie algebra is generated by invariant symmetric bilinear form or is trivial
depending on $A$ is symmetrizable or not. The result was conjectured by Moody\cite{Moody_78} and assumed by Kac\cite{Kac_852}. As applications we discuss the rational homotopy types of Kac-Moody groups and their flag manifolds.
\end{abstract}

Key words: Cartan matrix, Weyl group, polynomial invariants of Weyl group, Kac-Moody Group, Flag manifold.

MSC(2010): Primary 17C99, secondly 55N45

\section{Introduction}

Let $A=(a_{ij})$ be an $n\times n$ integer matrix satisfying

(1) For each $i,a_{ii}=2$;

(2) For $i\not=j,a_{ij}\leq 0$;

(3) If $a_{ij}=0$, then $a_{ji}=0$.

\noindent then $A$ is called a Cartan matrix.

Let $h$ be the real vector space spanned by $\Pi^{\vee}=\{\alpha^{\vee}_1,\alpha^{\vee}_2,\cdots,\alpha^{\vee}_n\}$, denote the dual basis of $\Pi^\vee$ in vector space $h^*$ by $\{\omega_1,\omega_2,\cdots,\omega_n\}$. That is $\omega_i(\alpha_j^{\vee})=\delta_{ij}$ for $1\leq i,j\leq n$. Let $\Pi=\{\alpha_1,\cdots,\alpha_n\}\subset h^*$ be given by $\langle\alpha^{\vee}_i,\alpha_j\rangle=a_{ij}$ for all $i,j$, then $\alpha_i=\sum\limits_{j=1}^n a_{ji}\omega_j$. Note that if the Cartan matrix $A$ is singular, then $\alpha_i,1\leq i\leq n$ is not a basis of $h^*$. $\Pi$ and $\Pi^\vee$ are called the root system and coroot system associated to Cartan matrix $A$. $\alpha_i, \alpha_i^\vee, \omega_i,1\leq i\leq n$ are called respectively the simple roots, simple coroots and fundamental dominant weights.

By the work of Kac\cite{Kac_68} and Moody\cite{Moody_68}, it is well known that for each Cartan matrix $A$, there is a Lie algebra $g(A)$ associated to
$A$ which is called Kac-Moody Lie algebra.


The Kac-Moody Lie algebra $g(A)$ is generated by $\alpha^\vee_i,e_i,f_i,1\leq i\leq n$ over $\mathbb{C}$,
with the defining relations:

(1) $[\alpha_i^{\vee},\alpha_j^{\vee}]=0$;

(2) $[e_i,f_j]=\delta_{ij}\alpha_i^{\vee}$;

(3) $[\alpha_i^{\vee},e_j]=a_{ij}e_j,[\alpha_i^{\vee},f_j]=-a_{ij}f_j$;

(4) $\mathrm{ad}(e_i)^{-a_{ij}+1}(e_j)=0$;

(5) $\mathrm{ad}(f_i)^{-a_{i j}+1}(f_j)=0$.

Kac and Peterson\cite{Kac_Peterson_83}\cite{Kac_Peterson_84}\cite{Kac_85} constructed the Kac-Moody group $G(A)$ with Lie algebra $g(A)$. In this paper for convenience we consider
the quotient Lie algebra  of $g(A)$ modulo its center $c(g(A))$ and the associated simply connected group $G(A)$ modulo $C(G(A))$. We still use the same symbols $g(A)$ and $G(A)$ and call
them the Kac-Moody Lie algebra and the Kac-Moody group.

Cartan matrices and their associated Kac-Moody Lie algebras, Kac-Moody groups are divided into three types.

(1) Finite type, if $A$ is positive definite. In this case, $G(A)$ is just the simply connected complex semisimple Lie group with Cartan matrix $A$.

(2) Affine type, if $A$ is positive semi-definite and has rank $n-1$.

(3) Indefinite type otherwise.

A Cartan matrix $A$ is called hyperbolic if all the proper principal sub-matrices of $A$ are of finite or affine type. A Cartan matrix $A$ is called symmetrizable if there exist an invertible diagonal matrix $D$ and a symmetric matrix $B$ such that $A=DB$. $g(A)$ is called a symmetrizable Kac-Moody Lie algebra if $A$ is symmetrizable.

The Weyl group $W(A)$ associated to Cartan matrix $A$ is the group generated by the Weyl reflections $\sigma_i:h^*\to h^*$ with respect to $\alpha_i^{\vee},1\leq i\leq n$, where $\sigma_i(\alpha)=\alpha-\langle \alpha,\alpha_i^\vee\rangle\alpha_i$. $W(A)$ has a Coxeter presentation
$$W(A)=<\sigma_1,\cdots,\sigma_n|\sigma^2_i=e,1\leq  i\leq n; (\sigma_i\sigma_j)^{m_{ij}}=e,1\leq i<j\leq n>.$$ where $m_{ij}=2,3,4,6 $ and $\infty$ as $a_{ij}a_{ji}=0,1,2,3 $ and $\geq 4$ respectively. The action of $\sigma_i$ on fundamental dominant weights is given by $\sigma_i(\omega_j)=\omega_j-(\omega_j,\alpha_i^{\vee})\alpha_i=\omega_j-\delta_{ji}\alpha_i$. For details see Kac\cite{Kac_82}, Humphureys\cite{Humphreys_90}.

The action of Weyl group $W(A)$ on $h^*$ induces an action of $W(A)$ on the polynomial ring $\mathbb{Q}[h^*]\cong \mathbb{Q}[\omega_1,\cdots,\omega_n]$. $f\in \mathbb{Q}[h^*]$ is called a $W(A)$ invariant polynomial if for each $\sigma\in W(A),\sigma(f)=f$. Since $W(A)$ is generated by $\sigma_i,1\leq i\leq n$, $f$ is a $W(A)$ invariant polynomial if and only if $\sigma_i(f)=f$ for $1\leq i\leq n$.
All the $W(A)$ invariant polynomials form a ring, called the ring of $W(A)$ polynomial invariants, denoted by $I(A)$.

The invariant theory of Weyl groups has been a significant topic since the 1950s. It has important applications in the homology of Lie groups and their classifying spaces. Motivated by that study, Chevalley showed that the ring of invariants of a finite Weyl group is a polynomial algebra. A comprehension study of the polynomial invariants was undertaken by Bourbaki, Solomon, Springer and Steinberg, etc.


In his paper \cite{Moody_78} Moody proved the following theorem.

\noindent {\bf Theorem}(Moody) Let A be an indecomposable and symmetrizable $n\times n$
Cartan matrix which associated invariant bilinear form $\psi$ is non-degenerate and of signature
$(n-1,1)$, then the ring of $W(A)$ polynomial invariants is $\mathbb{Q}[\psi]$.

In the same paper, Moody further said: \lq\lq
{\sl We conjecture that it is in fact true for all
Weyl groups arising from non-singular Cartan matrices of non-finite
type}\rq\rq

In \cite{Kac_852}, 
Kac also assumed that for an indecomposable and indefinite Cartan matrix
the ring of $W(A)$ polynomial invariants is $\mathbb{Q}[\psi]$ or trivial depending on $A$ is symmetrizable or not.

In this paper, we prove the following theorem.

\noindent{\bf Theorem: }Let $A$ be an indecomposable and indefinite Cartan matrix $A$. If $A$ is symmetrizable, then $I(A)=\mathbb{Q}[\psi]$; If $A$ is non symmetrizable, then $I(A)=\mathbb{Q}$.

The content of this paper is as below. In section 2, we discuss the general results about the polynomial invariants of Weyl group for a Kac-Moody group. In section 3 and 4 we consider the rank 2 case and the hyperbolic case respectively. The main theorem is proved in section 5. In section 6, we consider the applications of the theorem in determining the rational homotopy type of Kac-Moody groups and their flag manifolds.

\section{Rings of polynomial invariants of Weyl groups, general case}

In this section, we discuss some general properties of the ring of invariants of Weyl groups.

\noindent\begin{lemma} If a Cartan matrix $A$ is the direct sum of Cartan matrices $A_1,A_2$, then $I(A)\cong I(A_1)\otimes I(A_2)$.\end{lemma}

So we only consider indecomposable Cartan matrices.

\noindent\begin{lemma} Let $f\in I(A)$, $f=\sum\limits_{i=0}^l f_i$, and $f_i$ be the degree $i$ homogeneous component of $f$, then $f_i\in I(A)$.\end{lemma}

So $I(A)$ is a graded ring and $I(A)=\bigoplus\limits_{i=0}^{\infty} I^i(A)$. To determine the ring $I(A)$, we only need to consider homogeneous invariant polynomials.

\noindent\begin{lemma} For an indecomposable Cartan matrix $A$ of affine type or indefinite type, the orbit $\{\sigma(\omega)|\sigma\in W(A)\}$ of an element $\omega$ in Tits cone is an infinite set.\end{lemma}

\noindent{\bf Proof: }Since $\omega$ is in Tits cone, we can assume $\omega=\sum\limits_{i=1}^n \lambda_i \omega_i,\lambda_i\geq 0$. Let $S=\{1,2,\cdots,n\}, I=\{i\in S|\lambda_i=0\}$. Since $\{\sigma(\omega)|\sigma\in W(A)\}\cong W(A)/W_I(A)$ indexes the Schubert varieties of generalized flag manifold $F(A_I)=G(A)/G_I(A)$. The result is get from the fact that the number of Schubert varieties in $F(A_I)$ is infinite for affine and indefinite type.

\noindent\begin{coro}Let $f\in I(A)$ be a homogeneous invariant polynomial, $\omega\not=0$ is in the Tits cone. If $\omega|f$, then $f=0$ \label{w_n|f}.\end{coro}

\noindent{\bf Proof: }If $\omega|f$, then for any $\sigma\in W(A),\sigma(\omega)|\sigma(f)=f$. Since $\{\sigma(\omega)|\sigma\in W(A)\}$ is an infinite set, if $f\not=0$, this
contradicts to the condition that the degree of $f$ is finite.

\noindent\begin{lemma} For a Cartan matrix $A$, $I^1(A)=\{0\}$.\end{lemma}

\noindent {\bf Proof: }Suppose $f=\sum\limits_{i=1}^n \lambda_i \omega_i\in I^1(A)$, then for each $j$, $\sigma_j(f)=f-\lambda_j\alpha_j=f$. Since $\alpha_j\not=0$, $\lambda_j=0$. Therefore $f=0$.

\noindent\begin{lemma}Let $A$ be a Cartan matrix, $\lambda_{ij}=\lambda_{ji}$, then $f(\omega)=\sum\limits_{i,j=1}^n\lambda_{ij}\omega_i\omega_j \in I^2(A)$ is an $W(A)$ invariant polynomial if and only if $\displaystyle{\frac{\partial f}{\partial \omega_j}=\frac{1}{2}\frac{\partial^2 f}{\partial \omega_j^2}\alpha_j}$ for all $j$. That is $2\lambda_{ij}=a_{ij}\lambda_{jj}$ for all $i,j$.\label{relation of lambda and a}\end{lemma}

\noindent {\bf Proof: }If $f$ is an invariant polynomial, then for each $j$, $$\sigma_j(f)=f(\omega_1,\cdots,\omega_j-\alpha_j,\cdots,\omega_n)=f(\omega)-\frac{\partial f}{\partial \omega_j}\alpha_j+\frac{1}{2}\frac{\partial^2 f}{\partial \omega_j^2}\alpha_j^2=f.$$
It is equivalent to $\displaystyle{\frac{\partial f}{\partial \omega_j}=\frac{1}{2}\frac{\partial^2 f}{\partial \omega_j^2}\alpha_j}$. That is  $2\sum\limits_{i=1}^n\lambda_{ij}\omega_i=\lambda_{jj}\alpha_j=\lambda_{jj}\sum\limits_{i=1}^n a_{ij} \omega_i$, i.e. $2\lambda_{ij}=a_{ij}\lambda_{jj}$.

Lemma \ref{relation of lambda and a} can be generalized to

\noindent\begin{lemma}Let $A$ be a Cartan matrix, $f$ is a degree $l$ invariant polynomial if and only if for all $j$ $$\displaystyle{\frac{\partial f}{\partial \omega_j}-\frac{1}{2!}\frac{\partial^2 f}{\partial \omega_j^2}\alpha_j}+\cdots+(-1)^l \frac{1}{l!}\frac{\partial^l f}{\partial \omega_j^l}\alpha^{l-1}_j=0.\hspace{4cm} (1)$$\end{lemma}

\noindent \begin{lemma} An $n\times n$ Cartan matrix $A$ is symmetrizable, if and only if there exist non-zero $d_1,d_2,\cdots,d_n$ such that $a_{ij}d_j=a_{ji}d_i$ for all $i,j$. \label{a and d}\end{lemma}

\noindent{\bf Proof: }Suppose $A$ is symmetrizable, then there exist an invertible diagonal matrix $D=\mathrm{diag}(d_1,\cdots,d_n)$ and a symmetric matrix $B$, such that $A=DB$, that is $a_{ij}=d_i b_{ij}$ for all $i,j$. So $\displaystyle{\frac{a_{ij}}{d_i}=b_{ij}=b_{ji}=\frac{a_{ji}}{d_j}}$. It is equivalent to $a_{ij}d_j=a_{ji}d_i$.

If there exist non-zero $d_1,d_2,\cdots,d_n$ such that $a_{ij}d_j=a_{ji}d_i$ for all $i,j$, then let $D=\mathrm{diag}(d_1,\cdots,d_n)$ and $B=(b_{ij})_{n\times n}=\displaystyle{(\frac{a_{ij}}{d_i})_{n\times n}}$, then $A=DB$. Therefore $A$ is symmetrizable.

\noindent \begin{coro} If $A$ is an indecomposable Cartan matrix, then $\dim I^2(A)=1$ or $0$ depending on whether $A$ is symmetrizable or not. And if $A$ is symmetrizable, then $I^2(A)$ is spanned by an invariant bilinear form $\psi$ which is unique up to a constant.\end{coro}

\noindent{\bf Proof: }For an indecomposable Cartan matrix, $\dim I^2(A)>0$ or $=0$.

If $\dim I^2(A)>0$, choose $f(\omega)=\sum\limits_{i,j=1}^n\lambda_{ij}\omega_i\omega_j \in I^2(A)$, $f(\omega)\not=0$. By permutating the simple roots, we can assume that there exists an integer $k> 0$, such that $\lambda_{11},\cdots,\lambda_{kk}\not=0$, but $\lambda_{k+1,k+1},\cdots,\lambda_{nn}=0$. If $i\leq k$ and $j>k$, by Lemma \ref{relation of lambda and a}, $2\lambda_{ij}=a_{ij}\lambda_{jj}$, therefore $\lambda_{ij}=0$. By $0=2\lambda_{ij}=2\lambda_{ji}=a_{ji}\lambda_{ii}$, we get $a_{ji}=0$ for all $i\leq k,j>k$. Since $A$ is indecomposable $k$ must equals to $n$. Let $d_i=\lambda_{ii}$ for $1\leq i\leq n$, then $a_{ij}d_j=a_{ji}d_i$ for all $i,j$, this means that $A$ is symmetrizable.

If $\dim I^2(A)=0$, by similar discussion, we know $A$ is non symmetrizable.

Since $A$ is indecomposable, $\lambda_{ij}:\lambda_{jj}$ and $\lambda_{ii}:\lambda_{jj}$ for all $i,j$ are determined by $A$. Therefore if $\dim I^2(A)>0$, it must be $\dim I^2(A)=1$.

Below for an indecomposable and symmetrizable Cartan matrix $A$, we always fix a non-zero $\psi\in I^2(A)$.

\section{Rings of polynomial invariants of Weyl groups, $n=2$ case }
A $2\times 2$ Cartan matrix is of form $\left(
                                         \begin{array}{cc}
                                           2 & -a \\
                                           -b & 2 \\
                                         \end{array}
                                       \right)$. We denote it by $A_{a,b}$. $A_{a,b}$ is of finite type(affine type or indefinite type) if $ab< 4$($ab=4$ or $ab>4$).  The action of reflections $\sigma_1,\sigma_2\in W(A)$ on $h^*$ is given by
$$\sigma_1(\omega_1)=-\omega_1+b\omega_2, \sigma_1(\omega_2)=\omega_2, \sigma_2(\omega_1)=\omega_1, \sigma_2(\omega_2)=-\omega_2+a\omega_1.$$
\noindent\begin{lemma}The Weyl group $W_{a,b}$ of Cartan matrix $A_{a,b}$ is dihedral group
$D_{m}$, where $m=2,3,4,6 $ and $\infty$ as $ab=0,1,2,3 $ and $\geq 4$ respectively. If $A_{a,b}$ is of affine or indefinite type, then
the ring of polynomial invariants of Weyl group $W_{a,b}$ is $I(A_{a,b})=\mathbb{Q}[\psi]$.\label{n=2}\end{lemma}






\noindent {\bf Proof: }For the Cartan matrix $A_{a,b}$ of affine or indefinite type, $ab\not=0$. Since $A$ is indecomposable and symmetrizable, $\dim I^2(A)=1$ and is spanned by $\psi=a\omega_1^2-ab\omega_1\omega_2+b\omega_2^2 \in I^2(A)$.
Suppose $f(\omega)=\sum\limits_{i=0}^l \lambda_i \omega_1^i\omega_2^{l-i}$ is a degree $l$ homogeneous invariant polynomial of degree $l$, then

$\displaystyle{\sigma_2(f)=\sum\limits_{i=0}^l \lambda_i \omega_1^i(-\omega_2+a\omega_1)^{l-i}}$

$\displaystyle{=\sum\limits_{i=0}^l \sum\limits_{j=0}^{l-i} (-1)^j \lambda_i {{l-i}\choose{j}}\omega_1^i\omega_2^j(a\omega_1)^{l-i-j}}$

$\displaystyle{=\sum\limits_{j=0}^l \sum\limits_{i=0}^{l-j} (-1)^j \lambda_i {{l-i}\choose{j}} a^{l-i-j} \omega_1^{l-j}\omega_2^j} $

$\displaystyle{=\sum\limits_{j=0}^l [\sum\limits_{i=0}^{j} (-1)^{l-j} \lambda_i {{l-i}\choose{l-j}} a^{j-i}] \omega_1^{j}\omega_2^{l-j}} $

So $\sigma_2(f)=f$ is equivalent to $$ \lambda_j=\sum\limits_{i=0}^{j} (-1)^{l-j} \lambda_i {{l-i}\choose{l-j}} a^{j-i},0\leq i\leq l\hspace{5cm} (2)$$
\indent Let $j=0$, we get $\lambda_0=(-1)^n \lambda_0 $. So $\lambda_0=0$ or $n$ is even.

1. If $\lambda_0=0$, then $\omega_1|f$. By Corollary \ref{w_n|f}, $f=0$.

2. If $l$ is even, suppose $l=2m$. There exists a constant $\lambda$ such that $f-\lambda \psi^m$ is an invariant polynomial and $\omega_1|(f-\lambda \psi^m)$, hence
$f=\lambda \psi^m$.

This proves the lemma.
\section{Some results about hyperbolic Cartan matrices}
In his paper \cite{Moody_78} Moody proved that for each indecomposable and symmetrizable hyperbolic Cartan matrix $A$, the ring of polynomial invariants $I(A)=\mathbb{Q}[\psi]$, where $\psi$ is the invariant bilinear form. So in this section we only consider non symmetrizable Cartan matrices.

The indecomposable hyperbolic Cartan matrix exists only for $n\leq 10$ and their numbers are finite for $3\leq n\leq 10$. There are lists of hyperbolic Cartan matrices in Wan \cite{Wan_91} and Carbone\cite{Carbone_10}.

\noindent \begin{lemma} Let $A$ be an indecomposable and non symmetrizable hyperbolic Cartan matrix with $n\geq 4$, then $A$ satisfies

C1. The Dynkin diagram of $A$ forms a circle. That is $a_{ij}\not=0$ if and only if $|i-j|=0,1$ or $n-1$.

C2. If $|i-j|=1$ or $n-1$, then $a_{ij}=-1$ or $a_{ji}=-1$.\label{C1C2}\end{lemma}

The lemma is proved by direct checking in the lists.

\begin{remark}The lemma is not true for the case $n=3$. \end{remark}

\noindent \begin{lemma} Let $A$ be an indecomposable and non symmetrizable hyperbolic Cartan matrix with $n=3$, then $A$ contains a $2\times 2$ principal sub-matrix of affine type or all the $2\times 2$ principal sub-matrices of $A$ are of finite type. In the latter case, $A$ satisfies the conditions C1, C2.\end{lemma}





In \cite{Feingold_Nicolai_04}, Feingold and Nicolai proved the following theorem.

\noindent{\bf Theorem: }Let $g(A)$ be a Kac-Moody Lie algebra associated to symmetrizable Cartan matrix $A=(a_{ij})_{n\times n}$
which is generated by $\alpha_i^\vee,e_i, f_i, 1\leq i\leq n$,
$\beta_1,\cdots,\beta_m$ be a set of positive real roots of $g(A)$ such that
$\beta_i-\beta_j,1\leq i\not=j\leq m$ are not roots. Let $E_i,F_i$ be root vectors in the one dimensional
root spaces corresponding to the positive real roots $\beta_i$ and the negative
real roots $-\beta_i$ respectively, and let $H_i = [E_i, F_i]$. Then the Lie subalgebra of $g$
generated by $\{E_i, F_i,H_i | 1 \leq i \leq m\}$ is a regular Kac-Moody subalgebra with Cartan matrix
$B=(b_{ij})_{n\times n}=(\frac{2(\beta_{j},\beta_i)}{(\beta_i,\beta_i)})_{n\times n}$.



By using the ideas in the theorem of Feingolds and Nicolai, we can prove the following lemma. We assume $n+1=1,(n-1)+2=1$, etc.

\noindent\begin{lemma}
Let $A$ be an $n\times n$ Cartan matrix satisfying the conditions C1 and C2 in Lemma \ref{C1C2} and its simple roots are $\{\alpha_1,\alpha_2,\cdots,\alpha_n\}$, then $\beta_i=\alpha_{i+1}+\alpha_{i+2}, 1\leq i\leq n$ is a set of positive real roots of $g(A)$ and $\beta_i-\beta_j$ for $i\not=j$ are not roots. Let $\alpha_i^\vee,e_i,f_i,1\leq i\leq n$ be the generators of $g(A)$, $E_i=[e_{i+1},e_{i+2}],F_i=-[f_{i+1},f_{i+2}]$ and $H_i=[E_i,F_i]$, then
$H_i,E_i,F_i,1\leq i\leq n$ generate a full rank regular Kac-Moody subalgebra with simple root system $\{\beta_1,\beta_2,\cdots,\beta_n\}$ and Cartan matrix $B=(b_{ij})=(\beta_{j}(H_i))_{n\times n}$.\label{g(B)}
\end{lemma}

\noindent{\bf Proof: }For $\beta_i=\alpha_{i+1}+\alpha_{i+2}$, 

$H_i=[E_i,F_i]$

$=-[[e_{i+1},e_{i+2}],[f_{i+1},f_{i+2}]]$

$=-[[[e_{i+1},e_{i+2}],f_{i+1}],f_{i+2}]-[f_{i+1},[[e_{{i+1}},e_{i+2}],f_{i+2}]]$

$=[[[f_{i+1} ,e_{{i+1}}],e_{i+2}],f_{i+2}]+[[e_{{i+1}},[f_{i+1}, e_{i+2}]],f_{i+2}]$

$+[f_{i+1},[[f_{i+2}, e_{{i+1}}],e_{i+2}]]+[f_{i+1},[e_{{i+1}},[f_{i+2},e_{i+2}]]]$

$=-[[\alpha_{i+1}^\vee,e_{i+2}],f_{i+2}]-[f_{i+1},[e_{{i+1}},\alpha_{i+2}^\vee]]$

$=-(a_{{i+1},{i+2}}\alpha_{i+2}^\vee+\alpha_{{i+2},{i+1}}\alpha_{i+1}^\vee)$.

\noindent then

$[H_i,E_i]=-2(a_{i+1,i+2}+a_{i+2,i+1}+a_{i+1,i+2}a_{i+2,i+1})E_i$










Hence for each $1\leq i\leq n$, $$-2(a_{i+1,i+2}+a_{i+2,i+1}+a_{i+1,i+2}a_{i+2,i+1})=2(1-(a_{i+1,i+2}+1)(a_{i+2,i+1}+1))=2, $$ A routine checking shows that $\beta_{j}(H_i)\leq 0$. Therefore the matrix $B$ with $b_{ij}=\beta_{j}(H_i)$ is a Cartan matrix. Hence $H_i,E_i,F_i,1\leq i\leq n$ generate a Kac-Moody Lie algebra with Cartan matrix $B$.

We can't say Lemma 4.3 is the corollary of the theorem, since we don't know whether the theorem is true for non symmetrizable Cartan matrix $A$. So we must prove Lemma 4.3 by direct computation.

\begin{coro}Let $A$ be an $n\times n$ indecomposable and non-symmetrizable hyperbolic Cartan matrix, then in $g(A)$ there is a full rank indecomposable and non hyperbolic regular indefinite Kac-Moody subalgebra $g(B)$. \label{regular subalgebra}\end{coro}

This corollary is proved by using the Lemma \ref{g(B)} and checking one by one in the list of indecomposable, non symmetrizable hyperbolic Cartan matrices. We have composed computer program to do the checking. The computation result shows except for the number 131,132,133,137,139,141 hyperbolic Lie algebras in the list of Carbone\cite{Carbone_10}, all the subalgebras we constructed are non symmetrizable.

Below is a simple example.

\noindent \begin{example}
For hyperbolic Cartan matrix $A=\left(
                    \begin{array}{ccc}
                      2 & -1 & -1 \\
                      -1 & 2 & -1 \\
                      -2 & -1 & 2 \\
                    \end{array}
                  \right)$, we get a regular subalgebra $g(B)$ of $g(A)$ with simple roots $\beta_1=\beta_2+\beta_3, \beta_2=\beta_3+\beta_1, \beta_3=\beta_1+\beta_2$, and the Cartan matrix is
$B=\left(
                    \begin{array}{ccc}
                      2 & -2 & -2 \\
                      -3 & 2 & -1 \\
                      -1 & -1 & 2 \\
                    \end{array}
                  \right)$.
It is non symmetrizable and indefinite.
\end{example}












\section{Proof of the main theorem}
\subsection{Some preparing lemmas}
Let $A$ be an $n\times n$ Cartan matrix and $S=\{1,2,\cdots,n\}$. For $I\subset S$, let $A_I$ be the principal sub-matrix $(a_{ij})_{i,j\in I}$ corresponding to $I$. Then $A'=A_{S-\{n\}}$ is the upper-left $(n-1)\times (n-1)$ principal sub-matrix of $A$. Let $h'$ be the subspace of $h$ spanned by $\alpha^{\vee}_1,\cdots,\alpha_{n-1}^{\vee}$ and $h'^*$ the subspace of $h^*$ spanned by $\omega_1,\cdots,\omega_{n-1}$, then $h=h'\oplus \mathbb{R} \alpha^{\vee}_n$ and $h^*=h'^*\oplus \mathbb{R} \omega_n$. Let $\alpha_i\in {h}^*,1\leq i\leq n$ and $\alpha'_i\in {h'}^*,1\leq i\leq n-1$ be respectively the simple roots of Cartan matrices $A$ and $A'$, $\sigma_i,1\leq i\leq n$ and $\sigma'_i,1\leq i\leq n-1$ be respectively the Weyl reflections on ${h}^*$ and ${h'}^*$. For $1\leq i\not=j\leq n-1$, $\alpha_i=\alpha'_i+a_{ni}\omega_n$, $\sigma_i(\omega_j)=\sigma'_i(\omega_j)$, and $\sigma_i(\omega_i)=\omega_i-\alpha_i=\sigma'_i(\omega_i)-a_{ni}\omega_n$.

\noindent\begin{lemma}Let $\omega=(\omega_1,\cdots,\omega_{n-1},\omega_{n}),\omega'=(\omega_1,\cdots,\omega_{n-1})$. $f(\omega)$ is a degree $l$ invariant polynomial under the action of
$\sigma_1,\cdots,\sigma_{n-1}$ and
$f(\omega)=\sum\limits_{i=0}^l f_i(\omega')\omega_n^{l-i}$ with $f_i(\omega')$ degree $i$ homogeneous polynomial in $S(h'^*)$, then $f_l(\omega')$ is invariant under the action of $\sigma'_1,\cdots,\sigma'_{n-1}$.\end{lemma}

\noindent {\bf Proof: }For $k\not=n$,

$f(\omega)=\sigma_k(f(\omega))$

$=\sum\limits_{i=0}^l \sigma_k(f_i(\omega'))\omega_n^{l-i}$

$=\sum\limits_{i=0}^l f_i(\sigma_k(\omega_1),\sigma_k(\omega_2),\cdots,\sigma_k(\omega_k),\cdots,\sigma_k(\omega_{n-1}))\omega_n^{l-i}$

$=\sum\limits_{i=0}^l f_i(\sigma'_k(\omega_1),\sigma'_k(\omega_2),\cdots,\sigma'_k(\omega_k)-a_{nk}\omega_n,\cdots,\sigma'_k(\omega_{n-1}))\omega_n^{l-i}$,

\noindent set $\omega_n=0$, we get $f_l(\omega')=f_l(\sigma'_k(\omega'))=\sigma'_k(f_l(\omega'))$.

\noindent\begin{coro}If $f(\omega)=\sum\limits_{i=0}^l f_i(\omega')\omega_n^{l-i}\in I(A)$, then $f_l(\omega')\in I(A')$ \label{fl invatiant}.\end{coro}

\noindent\begin{lemma}If the degree $l$ polynomial $f(\omega)$ is invariant under the action of
$\sigma_1,\cdots,\sigma_{n-1}$ and
$f(\omega)=\sum\limits_{i=0}^l f_i(\omega)\omega_n^{l-i}$, then $$\sigma'_k(f_i(\omega'))=\sum\limits_{j=0}^{l-i} \frac{(-a_{nk})^j}{j!} \frac{\partial ^j f_{i+j}(\omega')}{(\partial \omega_k)^j}\hspace{4cm} (3)$$\end{lemma}
\noindent {\bf Proof: }For $k\not=n$,

$\sum\limits_{i=0}^l f_i(\omega)\omega_n^{l-i}$

$=f(\omega)=\sigma_k(f(\omega))$

$=\sum\limits_{i=0}^l \sigma_k(f_i(\omega))\omega_n^{l-i}$

$=\sum\limits_{i=0}^l f_i(\sigma_k(\omega_1),\sigma_k(\omega_2),\cdots,\sigma_k(\omega_k),\cdots,\sigma_k(\omega_{n-1}))\omega_n^{l-i}$

$=\sum\limits_{i=0}^l f_i(\sigma'_k(\omega_1),\sigma'_k(\omega_2),\cdots,\sigma'_k(\omega_k)-a_{nk}\omega_n,\cdots,\sigma'_k(\omega_{n-1}))\omega_n^{l-i}$

$=\sum\limits_{i=0}^l \sum\limits_{j=0}^i \frac{1}{j!} \frac{\partial ^j f_i}{(\partial \omega_k)^j}(\sigma'_k(\omega_1),\sigma'_k(\omega_2),\cdots,\sigma'_k(\omega_k),\cdots,\sigma'_k(\omega_{n-1}))(-a_{nk}\omega_n)^j \omega_n^{l-i}$

$=\sum\limits_{i=0}^l \sum\limits_{j=0}^i \frac{(-a_{nk})^j}{j!} \frac{\partial ^j f_i}{(\partial \omega_k)^j}(\sigma'_k(\omega_1),\sigma'_k(\omega_2),\cdots,\sigma'_k(\omega_k),\cdots,\sigma'_k(\omega_{n-1})) \omega_n^{l-i+j}$

$=\sum\limits_{i=0}^l [\sum\limits_{j=0}^{l-i} \frac{(-a_{nk})^j}{j!} \frac{\partial ^j f_{i+j}}{(\partial \omega_k)^j}(\sigma'_k(\omega'))]\omega_n^{l-i}$

By comparing the coefficients of $\omega_n^{l-i}$ in the two sides, we get
$$f_i(\omega')=\sum\limits_{j=0}^{l-i} \frac{(-a_{nk})^j}{j!} \frac{\partial ^j f_{i+j}}{(\partial \omega_k)^j}(\sigma'_k(\omega')). $$
Acting $\sigma'_k$ on both sides we prove the lemma.







\noindent\begin{lemma}Let $A$ be an indefinite $n\times n$ Cartan matrix and its upper-left $(n-1)\times (n-1)$ principal sub-matrix is $A'$. If both $A$ and $A'$ are indecomposable and symmetrizable, then the restriction of the invariant bilinear form $\psi\in I(A)$ to $h'$ gives an invariant bilinear form $\psi'\in I(A')$.\label{psi psi}\end{lemma}

The proof is obvious by checking $\psi|_{h'}\not=0$ and $\psi'=\psi|_{h'}$ is invariant under the action of $\sigma'_1,\cdots,\sigma'_{n-1}$.

\noindent\begin{lemma}Let $f$ be a $W(A)$ invariant polynomial and $f(\omega)=\sum\limits_{i=0}^l f_i(\omega')\omega_n^{l-i}$, then $$f_j(\omega')=\sum\limits_{i=0}^{j} (-1)^{l-i} f_i(\omega'){{l-i}\choose l-j} {\omega'_n}^{j-i}\hspace{1cm}(4)\label{omega omega'}$$ where $\omega'_n=\sum\limits_{j\not=n}a_{jn} \omega_j$.\end{lemma}

\noindent {\bf Proof: }
$\sigma_n(\omega_n)=\omega_n-\alpha_n=-\omega_n-\sum\limits_{j\not=n}a_{jn} \omega_j=-\omega_n-\omega'_n$.

$f=\sigma_n(f)=\sum\limits_{i=0}^lf_i(\omega')\sigma_n(\omega_n^{l-i})$

$=\sum\limits_{i=0}^l f_i(\omega') (-\omega_n-\sum\limits_{j\not=n}a_{jn} \omega_j)^{l-i}$

$=\sum\limits_{i=0}^l f_i(\omega') \sum\limits_{j=0}^{l-i} (-1)^{l-i} {{l-i}\choose j}\omega^j_n {\omega'_n}^{l-i-j}$

$=\sum\limits_{j=0}^l [\sum\limits_{i=0}^{l-j} (-1)^{l-i} f_i(\omega'){{l-i}\choose j} {\omega'_n}^{l-i-j}]\omega^j_n$

$=\sum\limits_{j=0}^l [\sum\limits_{i=0}^{j} (-1)^{l-i} f_i(\omega'){{l-i}\choose l-j} {\omega'_n}^{j-i}]\omega^{l-j}_n$

By comparing the coefficients of  $\omega^{l-j}_n$, we prove the lemma.

\noindent\begin{remark}In fact the Equation (3) and (4) are just corollaries of Equation (1) applied to $f(\omega)=\sum\limits_{i=0}^l f_i(\omega')\omega_n^{l-i}$.\end{remark}
\noindent\begin{lemma}Let $f$ be a degree $l$ $W(A)$ invariant polynomial and $f(\omega)=\sum\limits_{i=0}^l f_i(\omega')\omega_n^{l-i}$. If $l=2m$, then for each $i\leq m-1$, there exist constants $a^i_j(l),0\leq j\leq i$ depending on $l$ such that
$f_{2i+1}=\sum\limits_{j=0}^i a^i_j(l) f_{2(i-j)}{\omega'}_n^{2j+1}$.
If $l=2m+1$, then $f_0=0$, and for each $i\leq m$, there exist constants $b^i_j(l),1\leq j\leq i$ such that
$f_{2i}=\sum\limits_{j=1}^i b^i_j(l) f_{2(i-j)+1}{\omega'}_n^{2j-1}$.
And the coefficients $a^i_j(l)$ and $b^i_j(l)$ can be computed.
\end{lemma}

\noindent {\bf Proof: }Let $j=0$ in Equation (4), we get $f_0=(-1)^lf_0$. So there are two cases.

Case 1, $l$ is even.

Let $j=1$ in Equation (4), we get $f_1=-f_1+{{l} \choose 1}f_0 \omega'_n$. That is $$f_1=\frac{1}{2} {l \choose 1} f_0 \omega'_n.\hspace{4cm} (5)$$

For $j=2$, we get $f_2=f_2 -{l-1 \choose 1} f_1 \omega'_n+ {{l} \choose 2} f_0 {\omega'}^2_n$, which is equivalent to $f_1=\frac{1}{2} {l \choose 1} f_0$.

For $j=3$, we get $$f_3=-f_3 +{l-2 \choose 1} f_2 \omega'_n- {{l-1} \choose 2} f_1 {\omega'}^2_n+  {{l} \choose 3} f_0 {\omega'}^3_n.$$
Substituting Equation (5) in it, we get
$$f_3=\frac{1}{2} {l-2 \choose 1} f_2 \omega'_n- \frac{1}{4}{{l} \choose 3} f_0 {\omega'}^3_n.\hspace{4cm} (6)$$
\indent Continuing this procedure, we prove the lemma when $l$ is even.


Case 2, $l$ is odd.

In this case, we have $f_0=0$. And the proof is similar to the case 1 that $l$ is even.

\noindent \begin{coro}Let $f$ be a $W(A)$ invariant polynomial and $f(\omega)=\sum\limits_{i=0}^l f_i(\omega')\omega_n^{l-i}$, then $\omega'_n|f_{l-1}(\omega')$\label{coro f_l-1}.\end{coro}

Computation motivate us to make the following conjecture.

\noindent {\bf Conjecture: }If $l$ is even, then the Equation (4) for $j=2k$ can be derived from the set of equations for $j=0,1,2,\cdots, 2k-1$. If $l$ is odd, the Equation (4) for $j=2k-1$ can be derived from the set of equations for $j=0,1,2,\cdots, 2k-2$.

The conjecture is verified for $k\leq 3$.


\noindent\begin{lemma}Let $A$ be an $n\times n$ Cartan matrix, if $f(\omega),g(\omega)\in S(h^*)$ satisfy $\sigma_k (f(\omega))-f(\omega)=\sigma_k (g(\omega))-g(\omega)$ for each $1\leq k\leq n$,
then $f-g\in I(A)$.\label{difference}\end{lemma}

The proof is trivial.

\noindent\begin{lemma}Let $A$ be an indefinite $n\times n$ Cartan matrix, its upper-left $(n-1)\times (n-1)$ principal sub-matrix is $A'$. If the ring of $W(A')$ polynomial invariants $I(A')=\mathbb{Q}[\psi']$ and $l=2m$, then for each $W(A)$ invariant polynomial $f(\omega)=\sum\limits_{i=0}^l f_i(\omega')\omega_n^{l-i}$ of degree $l$, there exists a constant $k$, such that $f_l(\omega')=k\psi'^{m}$ and $f_{l-1}(\omega')=km\psi'^{m-1}\omega_n^*$, here $\omega_n^*=\sum\limits_{k\not=n} \lambda_{kk} a_{nk}\omega_k$. \label{w^*}\end{lemma}

\noindent {\bf Proof: }By Corollary \ref{fl invatiant}, $f_l(\omega')\in I(A')$. Since $I(A')=\mathbb{Q}[\psi']$, there exists $k$, $f_l(\omega')=k\psi'^m$. In Equation (3), let $j=l-1$,
we get for $1\leq k\leq n-1$,
$$\sigma'_k(f_{l-1}(\omega'))-f_{l-1}(\omega')=-a_{nk} \frac{\partial f_{l}}{\partial \omega_k}.$$
Let $g(\omega')=k m\psi'^{m-1}\omega_n^*$, then it is easy to check $\displaystyle{-a_{nk} \frac{\partial f_{l}}{\partial \omega_k}=\sigma'_k(g(\omega'))-g(\omega')}$, so $\sigma'_k(f_{l-1}(\omega'))-f_{l-1}(\omega')=\sigma'_k(g(\omega'))-g(\omega')$. Applying Lemma \ref{difference} to $f_{l-1}, g$ for Cartan matrix $A'$, we get $f_{l-1}(\omega')-g(\omega')\in I^{l-1}(A')$. But $I^{l-1}(A')=I^{2m-1}(A')=\{0\}$, hence $f_{l-1}(\omega')=g(\omega')$.

\subsection{Proof of three propositions}
\noindent \begin{prop}Let $A$ be an $n\times n$ indecomposable and indefinite Cartan matrix, its upper-left $(n-1)\times (n-1)$ principal sub-matrix is $A'$. If
$I(A')=\mathbb{Q}$, then
$I(A)=\mathbb{Q}$.\end{prop}

\noindent {\bf Proof: }Let $f$ be a $W(A)$ invariant polynomial and $f(\omega)=\sum\limits_{i=0}^l f_i(\omega')\omega_n^{l-i}$, then by Corollary \ref{fl invatiant}, $f_l(\omega')\in I(A')$, so $f_l(\omega')=0$.

For $i=l-1$, the Equation (3) is
$$\sigma'_k(f_{l-1}(\omega'))=f_{l-1}(\omega')-a_{nk}\frac{\partial f_{l}}{\partial \omega_k}(\omega').$$
Substituting $f_l(\omega')=0$ in above equation, we get $\sigma'_k(f_{l-1}(\omega'))=f_{l-1}(\omega')$, so $f_{l-1}(\omega')=0$.
Continuing this procedure, we show that $f_i(\omega')=0$ for all $i>0$ and $f_0$ is a constant. Hence $f(\omega)=f_0 \omega_n^l$. By Corollary \ref{w_n|f}, $f=0$.

\noindent \begin{prop}Let $A$ be a symmetrizable and indefinite $n\times n$ Cartan matrix, its upper-left $(n-1)\times (n-1)$ principal sub-matrix is $A'$. If $I(A')=\mathbb{Q}[\psi']$, then $I(A)=\mathbb{Q}[\psi]$.\end{prop}

\noindent {\bf Proof: }Let $f$ be a $W(A)$ invariant polynomial and $f(\omega)=\sum\limits_{i=0}^l f_i(\omega')\omega_n^{l-i}$, then by Corollary \ref{fl invatiant}, $f_l(\omega')\in I(A')$, so $f_l(\omega')=0$ or there exists $\lambda\not=0, f_l=\lambda \psi'^m$. If $f_l=0$, then $\omega_n|f$, so $f=0$. If $f_l=\lambda \psi'^m$, then by Lemma \ref{psi psi} we can assume $\psi|_{h'}=\psi'$, so $f-\lambda \psi^m$ is a $W(A)$ invariant polynomial and $\omega_n|(f-\lambda \psi^m)$. Hence $f=\lambda \psi^m$.

\noindent \begin{prop}Let $A$ be an $n\times n$ indecomposable and non symmetrizable Cartan matrix, its upper-left $(n-1)\times (n-1)$ principal sub-matrix $A'$ is symmetrizable. If $I(A')=\mathbb{Q}[\psi']$, then $I(A)=\mathbb{Q}$.\end{prop}

\noindent {\bf Proof: }Let $f$ be a $W(A)$ invariant polynomial, $f(\omega)=\sum\limits_{i=0}^l f_i(\omega')\omega_n^{l-i}$. Suppose $\psi'=\sum\limits_{i,j=1}^{n-1} \lambda_{ij}\omega_i\omega_j, \lambda_{ij}=\lambda_{ji}$.

If $l$ is even, suppose $l=2m$. We prove $f_{l}=0$ at first. Suppose $f_{l}\not=0$, then by Lemma \ref{w^*}, there exists $k\not=0$, $f_l=k\psi'^{m}$ and $f_{l-1}=km\psi'^{m-1}\omega^*_n$. By Corollary \ref{coro f_l-1}, $\omega'_n|f_{l-1}=km\psi'^{m-1}\omega^*_n$. So $\omega'_n|\psi'$ or $\omega'_n|\omega_n^*$. Since $\psi'$ is $W(A)$ invariant and $-\omega'_n$ is in the Tits cone, by Corollary \ref{w_n|f}, $\omega'_n|\psi'$ is impossible. Therefore $\omega'_n|\omega_n^*$. Because $A$ is indecomposable, both $\omega'_n$ and $\omega_n^*$ are not $0$. Therefore there exists a constant $d_n\not=0$ such that $\omega_n^*=d_n \omega'_n$. But $\omega'_n=\sum\limits_{j\not=n}a_{jn} \omega_j$, $\omega_n^*=\sum\limits_{j\not=n} \lambda_{jj}a_{nj}\omega_j$, hence $a_{jn}d_n= a_{nj}\lambda_{jj}$. Let $d_i=\lambda_{ii},1\leq i\leq n-1$, since $A'$ is symmetrizable, by Lemma \ref{a and d} we know $a_{ij}d_j=a_{ji}d_i$ for all $i,j\leq n-1$. Combining with $a_{jn}d_n= a_{nj}\lambda_{jj}$, we get $a_{ij}d_j=a_{ji}d_i$ for all $i,j\leq n$. This shows $A$ is symmetrizable. But it is impossible. So $f_{l}=0$.

If $l$ is odd, then $f_l\in I^l(A')$ also implies $f_l=0$.

If $f_l=0$, then the remaining procedure of the proof is similar to the proof of Proposition 5.2.
\subsection{Proof of the main theorem}
To prove the main theorem we need the following lemma.

\noindent\begin{lemma}Let $A$ be a non-hyperbolic, indecomposable and indefinite Cartan matrix, then there exists an integer $k$, $1\leq k\leq n$ such that $A_{S-\{k\}}$ is an indecomposable and indefinite Cartan matrix.\label{find decomposable}\end{lemma}

\noindent {\bf Proof: }Since Cartan matrix $A$ is non-hyperbolic, there exists an integer $k$, $1\leq k\leq n$ such that $A_{S-\{k\}}$ is indefinite. If $A_{S-\{k\}}$ is indecomposable, the lemma is proved. If $A_{S-\{k\}}$ is decomposable, then the Dynkin diagram of $A_{S-\{k\}}$ is split to $r$ connected sub-diagram $\Gamma_1,\cdots,\Gamma_r$ with $r>1$ and there is a $s_0, 1\leq s_0\leq r$, such that the principal sub-matrix corresponding to $\Gamma_{s_0}$ is indefinite. Since $A$ is indecomposable, the simple root $\alpha_k$ is connected to all $\Gamma_s,1\leq s\leq r$.

By using the fact that the Dynkin diagrams of finite and affine types Cartan matrices are trees except for the case $A_n^{(1)}$, we can choose a suitable simple root $\alpha_{k'},k'\not=k$ and $\alpha_{k'}$ is not a vertex of Dynkin sub-diagram $\Gamma_{S_0}$, such that $A_{S-\{k'\}}$ is indecomposable and indefinite.

Now we can prove the main theorem.

\noindent{\bf Theorem: }Let $A$ be an indecomposable and indefinite Cartan matrix $A$. If $A$ is symmetrizable, then $I(A)=\mathbb{Q}[\psi]$; If $A$ is non symmetrizable, then $I(A)=\mathbb{Q}$.

\noindent {\bf Proof: }We prove this theorem by induction on $n$. For $n=2$, this is Lemma 3.1.

Suppose this theorem is true for all $(n-1)\times (n-1)$ indecomposable and indefinite Cartan matrices.

For an $n\times n$ indecomposable and indefinite Cartan matrices $A$, if $A$ is not hyperbolic, then by Lemma \ref{find decomposable}, we can find an $(n-1)\times (n-1)$ principal sub-matrix $A'$ which is both indecomposable and indefinite. Without loss of generality we can assume $A'$ is the upper-left $(n-1)\times (n-1)$ principal sub-matrix.


Then by considering the symmetrizability of $A'$ and $A$, there are three cases:

1. Both $A'$ and $A$ are non symmetrizable.

2. Both $A'$ and $A$ are symmetrizable.

3. $A'$ is symmetrizable and $A$ is non symmetrizable.

The proof for these three cases are dealt with by combining the induction assumption and Proposition 5.1, Proposition 5.2 and Proposition 5.3 respectively.

So except for the case $A$ is hyperbolic, we prove the theorem. For the hyperbolic case, if $A$ is symmetrizable, the proof is given in Moody\cite{Moody_78}; if $A$ is non symmetrizable, it is proved in the following proposition.

\noindent\begin{prop}
For an $n\times n$ indecomposable, non symmetrizable hyperbolic Cartan matrix $A$, $I(A)=\mathbb{Q}$.
\end{prop}

\noindent {\bf Proof: }For Cartan matrix $A$ with $n\geq 4$, by Lemma \ref{regular subalgebra} we can find an $n\times n$ indecomposable, non hyperbolic and indefinite Cartan matrix $B$ such that the root system associated to $B$ is a sub-root system of root system associated to $A$, and the Weyl group $W(B)$ is a subgroup of $W(A)$. Therefore $I(A)\subset I(B)$.

If $B$ is non symmetrizable, then by combining Lemma \ref{find decomposable}, Proposition 5.1 or 5.3 and the same induction procedure, we can prove $I(B)=\mathbb{Q}$. Hence $I(A)=\mathbb{Q}$.

If $B$ is symmetrizable, then by combining Lemma \ref{find decomposable} and Proposition 5.2, we prove $I(B)=\mathbb{Q}[\psi_B]$. To prove $I(A)=\mathbb{Q}$, it is sufficient to show $\psi_B^m,m\geq 1$ are not $W(A)$ invariants.

Suppose $\psi_B^m$ is a $W(A)$ invariant polynomial. If $m$ is odd, we get $\psi_B=(\psi_B^m)^{\frac{1}{m}}$ is $W(A)$ invariant. If $m$ is even, similarly we get for each $\sigma\in W(A)$, $\sigma(\psi_B)=\psi_B$ or $-\psi_B$. But $\sigma(\psi_B)=-\psi_B$ is impossible(A symmetric bilinear form $\psi=\sum\limits_{i,j=1}^n \lambda_{ij}\omega_i \omega_j$ with all the $\lambda_{ii},1\leq i\leq n$ having the same sign can't be transformed to $-\psi$ by a linear transformation). So we get $\sigma(\psi_B)=\psi_B$. Therefore $\psi_B$ is a $W(A)$ invariant polynomial. Since $A$ is non symmetrizable, this is impossible. Hence $I(A)=\mathbb{Q}$.

For the $n=3$ case, there are two possibilities. If $A$ contains a $2\times 2$ principal sub-matrix $A'$ of affine type, then by combining Lemma \ref{n=2} and Proposition 5.3, we show $I(A)=\mathbb{Q}$. If all the $2\times 2$ principal sub-matrices of $A$ are of finite type, then $A$ satisfies the conditions C1, C2. So we can find an indecomposable, non hyperbolic and indefinite Cartan matrix $B$ such that $g(B)$ is a regular subalgebra of $g(A)$. By the similar method for $n\geq 4$, we can also prove $I(A)=\mathbb{Q}$. This proves the proposition.

Thus the theorem is proved.








\section{Applications to rational homotopy types of Kac-Moody groups and their flag manifolds of indefinite type}


For the Kac-Moody Lie algebra $g(A)$, there is the Cartan decomposition $g(A)=h\oplus \sum\limits_{\alpha\in \Delta} g_{\alpha}$, where $h$ is the Cartan sub-algebra and $\Delta$ is the
root system of $g(A)$. Let $b=h\oplus \sum\limits_{\alpha\in \Delta^+} g_{\alpha}$ be the Borel sub-algebra, then $b$ corresponds to a Borel subgroup $B(A)$ in the Kac-Moody group $G(A)$. The
homogeneous space $F(A)=G(A)/B(A)$ is called the flag manifold of $G(A)$. By Kumar\cite{Kumar_02}, $F(A)$ is an ind-variety.

The cohomologies of Kac-Moody groups and their flag manifolds of finite and affine types are extensively studied. For reference see Pontrjagin\cite{Pontryagin_35}, Hopf\cite{Hopf_41}, Borel\cite{Borel_53_1}\cite{Borel_53}\cite{Borel_54}, Bott and Samelson\cite{Bott_Samelson_55}, Bott\cite{Bott_56}, Milnor and Moore\cite{Milnor_Moore_65} and Chevalley\cite{Chevalley_94}. But for indefinite type, little is known.

The rational cohomology rings of Kac-Moody groups and their flag manifolds are also considered in Kac\cite{Kac_852} and Kumar\cite{Kumar_85}. The essentially new part of our work is that we studied the properties of $P_A(q)$ and derive the explicit formula for $i_k$. For details see \cite{Jin_10}\cite{Jin_Zhao_11} and \cite{Zhao_Jin_12}.

For a Kac-Moody group $G(A)$, $H^*(G(A))$ is a locally finite free graded commutative algebra over $\mathbb{Q}$. Let the odd dimensional free generators of $H^*(G(A))$ be $y_1,\cdots,y_l$, and the even dimensional free generators of $H^*(G(A))$ be $z_1,\cdots,z_k,\cdots$. By Kac\cite{Kac_852}, Kichiloo\cite{Kitchloo_98} $l<n$. Denote the number of degree $k$ generators of $H^*(G(A))$ by $i_k$, then the Poincar\'{e} series of $G(A)$ is
$$
 P_G(q)=\prod\limits_{k=1}^{\infty}\frac{ (1-q^{2k-1})^{i_{2k-1}}}{(1-q^{2k})^{i_{2k}}}.
$$

The Poincar\'{e} series $P_G(q)$ determines the isometry type of the cohomology ring $H^*(G(A))$ and the rational homotopy type of $G(A)$.

Let $BB(A)$ be the classifying space of Borel subgroup $B(A)$ and $j: F(A)\to BB(A)$ be the classifying map of principal $B(A)$-bundle
$\pi:G(A)\to F(A)$. Denote the cohomology generators of $H^*(BB(A))$ by $\omega_1,\cdots,\omega_n,\deg \omega_i=2$.  A routine computation on the Leray-Serre spectral sequences of the fibration $G(A)\stackrel{\pi}{\longrightarrow} F(A) \stackrel{j}{\longrightarrow} BB(A)$ shows $$H^*(F(A))\cong E_3^{*,*}\cong \mathbb{Q}[\omega_1,\cdots,\omega_n]/<f_j,1\leq j\leq l>\otimes \mathbb{Q}[z_1,\cdots,z_k,\cdots].$$
where $f_j$'s correspond to the differential of $y_j$'s and they are the generators of the ring $I(A)$ of $W(A)$ polynomial invariants.

By previous work of the authors\cite{Zhao_Jin_12}, there is the following theorem.

\noindent\begin{theorem} Let $P_A(q)$ be the Poincar\'{e} series of flag manifold $F(A)$, then the sequence $i_2-i_1,i_4-i_3,\cdots,i_{2k}-i_{2k-1},\cdots$ can be derived from $P_A(q)$. In fact we can recover $P_A(q)$ from the sequence $i_2-i_1,i_4-i_3,\cdots,i_{2k}-i_{2k-1},\cdots$. \end{theorem}

But to determine the rational homotopy type of $G(A)$, we need to determine the sequence $i_1,i_2,\cdots,i_k,\cdots$. So except for the Poincar\'{e} series $P_A(q)$, we need more ingredients.
Note the number of generators of $I(A)$ of degree $k$ is just the integer $i_{2k-1}$. So if we can determine all the degrees of the generators in $I(A)$, then we can determine the sequence $i_1,i_3,\cdots,i_{2k-1},\cdots $. And the main theorem of this paper fills the gap. Now we have

\noindent \begin{theorem}For an indecomposable and indefinite Cartan matrix $A$, $i_{2k-1}=0$ for all $k>0$ except for $k=2$. And for $k=2$, if $A$ is symmetrizable, $i_3=1$; If $A$ is non symmetrizable, $i_3=0$. \end{theorem}


Set $\epsilon(A)=1$ or $0$ depending on $A$ is symmetrizable or not as in \cite{Kac_852}, then we get

\noindent \begin{theorem}The sequence $i_1,i_2,i_3,\cdots,i_{k},\cdots$ is determined from the Poincar\'{e} series $P_A(q)$ and $\epsilon(A)$.\end{theorem}

\noindent \begin{theorem}For an indecomposable and indefinite Cartan matrix $A$, the rational homotopy types of $G(A)$ are determined by the Poincar\'{e} series $P_A(q)$ and $\epsilon(A)$.\end{theorem}

Kumar\cite{Kumar_85} proved that for Kac-Moody Lie algebra $g(A)$, the Lie algebra cohomology $H^*(g(A),\mathbb{C})\cong H^*(G(A))\otimes \mathbb{C}$, so we also computed $H^*(g(A),\mathbb{C})$.
For a Kac-Moody group $G(A)$, $i_1=i_2=0$. And we have

\noindent \begin{coro}For an indecomposable and non symmetrizable indefinite Cartan matrix $A$, $G(A)$ is a 3-connected spaces.\end{coro}

\noindent \begin{coro}The dimension of the odd rational homotopy group $\pi_{odd} (G(A))$ of an indefinite Kac-Moody group $G(A)$ is $1$ or $0$ depending on $A$ is symmetrizable or not.\end{coro}

\noindent\begin{theorem}For an indecomposable and indefinite Cartan matrix $A$, if $A$ is symmetrizable, then
$$H^*(G(A))\cong \Lambda_{\mathbb{Q}}(y_3)\otimes \mathbb{Q}[z_1,\cdots,z_k,\cdots]$$ and $$H^*(F(A))\cong  \mathbb{Q}[\omega_1,\cdots,\omega_n]/<\psi>\otimes \mathbb{Q}[z_1,\cdots,z_k,\cdots]. $$
If $A$ is non symmetrizable, then
$$H^*(G(A))\cong \mathbb{Q}[z_1,\cdots,z_k,\cdots]$$ and $$H^*(F(A))\cong  \mathbb{Q}[\omega_1,\cdots,\omega_n]\otimes \mathbb{Q}[z_1,\cdots,z_k,\cdots]. $$
where $\deg z_k\geq 4$ are even for all $k$ and their degrees can be determined from the Poincar\'{e} series $P_A(q)$ and $\epsilon(A)$.\end{theorem}

Note the Poincar\'{e} series $P_A(q)$ can be computed easily by an inductive procedure. See \cite{Jin_10}\cite{Jin_Zhao_11} for details. So in principle the computation of rational homotopy type is solved for all the indecomposable and indefinite Kac-Moody groups, despite symmetrizable or not.

Since Kac-Moody groups and their flag manifolds are products of indecomposable Kac-Moody groups and indecomposable Kac-Moody flag manifolds, by combining with the known results for finite and affine types, we have determined the rational homotopy types of all the Kac-Moody groups and their flag manifolds. Since $G(A)$ and $F(A)$ are rational formal, see Sullivan\cite{Sullivan_77} and Kumar\cite{Kumar_02}, the rational homotopy groups and the rational minimal model of the corresponding Kac-Moody group $G(A)$ and its flag manifold $F(A)$ can be directly computed from Theorem 6.5.

\noindent\begin{theorem}For an $n\times n$ indecomposable and indefinite Cartan matrix $A$ satisfying $a_{ij}a_{ji}\geq 4$, the rational homotopy type of $G(A)$ is determined by $\epsilon(A)$.\label{>4}\end{theorem}

Since there are a large number of Cartan matrices satisfying the condition of Theorem \ref{>4}, this assertion seems to be very crazy. But the proof is very simple. It is derived from $P_A(q)=\frac{1+q}{1-(n-1)q}$.

It is deserved to mention that for $3\times 3$ non symmetrizable Cartan matrix $A$ with $a_{ij}a_{ji}\geq 4$ for all $i,j$, the Kac-Moody group $G(A)$ is a $5$ connected space. 

For an indecomposable and symmetrizable Cartan matrix $A$, let $p,q,r$ be the dimensions of positive, negative and zero vector subspaces of the invariant bilinear form $\psi$, set $\tau(A)=(p,q,r)$.

\noindent \begin{theorem} For an indecomposable and indefinite Cartan matrix $A$, if $A$ is symmetrizable, then the cohomology ring $H^*(F(A),\mathbb{C})$ is determined by $P_A(q)$ and $\tau(A)$.
If $g(A)$ is non symmetrizable, then the cohomology ring $H^*(F(A),\mathbb{C})$ is determined by $P_A(q)$. \end{theorem}

This is got from the Theorem 6.5 and classification of real quadratic forms.









\noindent {\bf Acknowledgements: } We would like to thank Professor Feingold and Nicolai for clearing some facts about the Theorem 3.1 in their paper.

\end{document}